\documentclass[11pt]{article}

\usepackage{amsmath}
\usepackage{amsfonts}
\usepackage[portuges]{babel}
\usepackage{amssymb}
\usepackage{latexsym}
\usepackage{enumerate}
\usepackage[dvips]{graphicx}

\def\rain{\to \infty}

\def\N{{\rm I\kern-.20em N}}
\def\R{{\rm I\kern-.20em R}}
\def\indi{{1\kern-.20em\rm I}}

\def\bkR{{\rm I\kern-.17em R}}
\def\bkN{{\rm I\kern-.20em N}}

\newtheorem{proposition}{Proposition}[section]

\newtheorem{definition}{Definition}[section]
\newtheorem{example}{Example}[section]

\newcommand{\bdem} {\begin{proof}}
\newcommand {\edem}{\hfill $\square$ \end {proof}}

\begin{document}
  \title{Dependence of multivariate extremes}
\author{C. Viseu \footnote{Instituto Superior de Contabilidade e
Administra\c c\~ ao de Coimbra, Coimbra,
Portugal.\newline
E-mail:cviseu@iscac.pt}\\  Instituto Superior de Contabilidade\\
e Administra\c c\~ ao de Coimbra\\ Portugal \and L. Pereira, A.P.
Martins, H. Ferreira \footnote{Departamento de
Matem\'atica, Universidade da Beira Interior, 6200 Covilh\~a,
Portugal.
E-mail: lpereira@ubi.pt; amartins@ubi.pt;
helena.ferreira@ubi.pt}\\Departamento de Matem\'atica \\
Universidade da Beira Interior\\ Portugal}
\date{}
\maketitle

\noindent {\bf Abstract:} We give necessary and sufficient conditions
for two sub-vectors of a random vector with a multivariate extreme value distribution,
corresponding to the limit distribution of the maximum of a multidimensional stationary sequence with extremal index, to be independent or totally dependent.
Those conditions involve first relations between the multivariate extremal
indexes of the sequences and secondly a coefficient that measure the strength of dependence between both sub-vectors.
The main results are illustrated with an auto-regressive sequence and a 3-dependent sequence.

\section{Introduction}
\label{sec:1}

Multivariate Extreme Value Analysis is frequently applied in the context of modeling environmental data, for which the phenomenon of dependence is often intrinsic. This paper focuses on the characterization of total dependence and of independence of two multivariate extreme value distributions.

Let ${\bf{X}}=\{{\bf{X}}_{n}^{(d)}=(X_{n,1},\ldots
,X_{n,d}\}_{n\geq 1}$ be a
$d-$dimensional stationary sequence with common distribution
function (d.f.) $Q({\bf{x}}^{(d)})=Q(x_{1},\ldots ,x_{d}),$
$\mathbf{x}^{(d)}\in
\mathbb{R}^{d},$ and  ${\bf{M}}_{n}=({M}_{n,1},\ldots ,{M}_{n,d}%
)$ the vector of pointwise maxima, where
${M}_{n,i}$ is the maximum of $i$-th component of $\bf X$.
Denote $\widehat{\mathbf{M}}_{n}=(\widehat{M}_{n,1}, \ldots
,\widehat{M}_{n,d})$ the corresponding vector of pointwise maxima of the
associated $d$-dimensional sequence, $\widehat{\mathbf{X}}=\{\widehat{\mathbf{X}}_{n}^{(d)}\}_{n\geq
1}$, of independent and identically distributed
(i.i.d.) random vectors having the same distribution function
$Q.$

In this multivariate setting operations among vectors are defined componentwise, that is
for each $d>1$ and $\mathbf{a}^{(d)}, \mathbf{b}^{(d)}\in\mathbb{R}^d$, $\mathbf{a}^{(d)}\leq\mathbf{b}^{(d)}$, if and only if $a_{j}\leq b_{j}$, for all $j=1,2,\ldots,d$.

If there exist sequences $\{(a_{n,1}>0, \ldots,a_{n,d}>0)\} _{n\geq 1}$ and $\{
(b_{n,1},\ldots ,b_{n,d})\}_{n\geq 1},$ such
that for
$\mathbf{u}(\mathbf{x}^{(d)})=\{\mathbf{u}_{n}(\mathbf{x}^{(d)})=(
a_{n,1}x_{1}+b_{n,1},\ldots,a_{n,d}x_{d}+b_{n,d}) \} _{n\geq
1},$
$$P\left( \widehat{\mathbf{M}}_n\leq \mathbf{u}_{n}(\mathbf{x}
^{(d)})\right) =P\left( \underset{j=1}{\overset{d}{\bigcap
}}\left\{ \widehat{M}_{n,j}\leq a_{n,j}x_{j}+b_{n,j}\right\} \right)
\xrightarrow [n\rain]{}G(\mathbf{x}^{(d)}),\ \mathbf{x}^{(d)}\in
\mathbb{R}^{d},$$ where $G$ is a d.f. with non-degenerate margins,
then $Q$ is said to be in the max-domain of
attraction of $ G$ ($Q\in D(G)$) and $G$ is said to be a Multivariate Extreme
Value (MEV) distribution function.

We will assume, without loss of generality, that the univariate
marginal distributions of $G$ are equal to $F$.

It is well known that the relationship between the d.f. $G(\mathbf{x}^{(d)})$, $\mathbf{x}^{(d)}\in
\mathbb{R}^{d},$ and its marginal distributions $F(x_j),j=1,\ldots
,d,$ can be characterized by its copula function, $D_G,$ which exhibits a
number of interesting properties, namely its stability equation:
\begin{equation}
D_{G}^{t}(y_{1},\ldots ,y_{d})=D_{G}(y_{1}^{t},\ldots ,y_{d}^{t}),\text{ }\forall t>0\text{  and  }\left(
y_{1},\ldots ,y_{d}\right) \in \left[ 0,1\right] ^{d}.  \label{eq1}
\end{equation}

If the stationary sequence ${\bf {X}}$ satisfies some mixing conditions, $D(\mathbf{u}_n(\mathbf{x}^{(d)}))$ of
Hsing (\cite{Hsing}) or
 $\Delta (\mathbf{u}_n(\mathbf{x}^{(d)}))$ of Nandagopalan (\cite{Nandagopalan}), and
\begin{equation*}
P\left( \mathbf{M}_n\leq \mathbf{u}_n(\mathbf{x}^{(d)})\right)
\xrightarrow [n\rain]{}H(\mathbf{x}^{(d)}),\ \mathbf{x}^{(d)}\in
\mathbb{R}^{d},
\end{equation*}
where $H$ is a d.f. with non-degenerate components, then $H$ is also a MEV d.f..
The MEV d.f $H$ and $G$ can be related through
the multivariate extremal index function, $\theta({\boldsymbol{\tau}}^{(d)})=\theta\left(\tau_1,\ldots,\tau_d\right)$ introduced by Nandagopalan (\cite{Nandagopalan}), which is a measure of clustering among the extreme values of a multivariate stationary sequence.

\begin{definition}
A $d-$dimensional stationary sequence $\mathbf{X}$,
is said to have multivariate extremal index $\theta^{\mathbf{X}}
({\boldsymbol{\tau}}^{(d)}) \in \left[
0,1\right] ,$ if for each $
\boldsymbol{\tau}^{(d)}=\left( \tau_1,\ldots ,\tau_d\right) \in \mathbb{R}_+^{d}$ there exists
$\mathbf{u}_{n}^{(
\boldsymbol{\tau}^{(d)})}=(u_{n,1}^{(\tau_1)},\ldots
,u_{n,d}^{(\tau_{d})}),$ $n\geq 1,$ satisfying
\begin{equation*}
nP(X_{1,j}>u_{n,j}^{(\tau_j)})\xrightarrow [n\rain]{}\tau_j, \ j=1,2,\ldots,d, \quad
P\left(\widehat{\mathbf{M}}_{n}\leq \mathbf{u}_{n}^{(
\boldsymbol{\tau}^{(d)})}\right) \xrightarrow [n\rain]{}G ({\boldsymbol{\tau}}^{(d)}) \quad \textrm{and}
\end{equation*}
\begin{equation*}
P\left(\mathbf{M}_{n}\leq \mathbf{u}_{n}^{({\boldsymbol{\tau}}^{(d)})}\right)
\xrightarrow [n\rain]{}G^{\theta^{\mathbf {X}} ({\boldsymbol{\tau}}^{(d)})} ({\boldsymbol{\tau}}^{(d)}).
\end{equation*}
\end{definition}

As in one dimension, the extremal index is a key parameter when relating the properties of extreme values of a stationary sequence to those of independent random vectors from the same $d-$dimensional marginal distribution. However, unlike the one dimensional case, it is not a constant for the whole process, but instead depends on the vector ${\boldsymbol{\tau}}^{(d)}$.

It is now clear that the existence of $\theta^{\mathbf{X}}
(\boldsymbol{\tau}^{(d)})$ allows us to write
\begin{equation*}
H({\bf x}^{(d)})=G^{\theta^{\bf X}({\boldsymbol{\tau}}^{(d)})}({\bf x}^{(d)}) \quad \textrm{with} \quad \tau_j\equiv\tau_j(x_j)=-\log F(x_j), \ j=1,2,\ldots,d.
\end{equation*}

Taking $d=p+q$, it follows, as a consequence of the definition of multivariate extremal index, that the sequences ${\bf{X}}^{(p)}=\{ \textbf{X}_n^{(p)}=(X_{n,1},\ldots ,X_{n,p})\}_{n\geq
1}$ and $\mathbf{X}^{(q)}=\{\mathbf{X}_n^{(q)}=(X_{n,p+1},\ldots
,X_{n,p+q})\}_{n\geq 1}$ have, respectively, extremal indexes
$$\theta^{\mathbf{X}^{(p)}}(\boldsymbol{\tau}^{(p)})=\lim_{{\tau_j\rightarrow 0^+}\atop{j=p+1,\ldots ,p+q}}\theta^{\mathbf{X}}(\boldsymbol{\tau}^{(p+q)})\quad \textrm{and}\quad \theta^{\mathbf{X}^{(q)}}(\boldsymbol{\tau}^{(q)})=\lim_{{\tau_j\rightarrow 0^{+}}\atop{j=1,\ldots ,p}}\theta^{\mathbf{X}}(\boldsymbol{\tau}^{(p+q)}).$$

In the notation of the extremal index we shall omit the sequence, whenever it is clear by
the context and the argument of the function.

Hereinafter, let ${\bf{Y}}=(Y_1,\ldots,Y_{p+q})$ and $\widehat{\bf {Y}}=(\widehat{Y}_1,\ldots,\widehat{Y}_{p+q})$ be, respectively,\linebreak  two random vectors
with distribution functions $G^{\theta^{\mathbf{X}}\left({\boldsymbol{\tau}}^{(p+q)}\right)}$ and $G$, where
$\mathbf{Y}^{(p)}=(Y_1,\ldots,Y_p)$ and $\mathbf{Y}^{(q)}=(Y_{p+1},\ldots,Y_{p+q})$ denote two sub-vectors of ${\bf{Y}}$ and
$\widehat{\mathbf{Y}}^{(p)}=(\widehat{Y}_1,\ldots,\widehat{Y}_p)$ and $ \widehat{\mathbf{Y}}^{(q)}=(\widehat{Y}_{p+1},\ldots,\widehat{Y}_{p+q})$
two sub-vectors of $\widehat{\mathbf {Y}}$.

In section 2 we discuss conditions under which $\mathbf{Y}^{(p)}$ and $\mathbf {Y}^{(q)}$ are independent or totally dependent. These conditions are established first by relations between the extremal indexes  $\theta^{\mathbf{X}}(\boldsymbol{\tau}^{(p+q)})$, $\theta^{\mathbf{X}^{(p)}}(\boldsymbol{\tau}^{(p)})$ and $\theta^{\mathbf{X}^{(q)}}(\boldsymbol{\tau}^{(q)})$ and secondly by a coefficient that measures the strength of dependence between $\mathbf{Y}^{(p)}$ and $\mathbf{Y}^{(q)}$.

The main results are illustrated in section 3 with an auto-regressive sequence and a 3-dependent sequence.

\section{Main results}
\label{sec:2}

If the d.f. $Q$ belongs to the domain of attraction of a MEV
distribution, $G$, and ${\bf{X}}$ has extremal index $\theta({\boldsymbol{\tau}}^{(p+q)}),$ ${\boldsymbol{\tau}}^{(p+q)}=(\tau_1,\dots,\tau_{p+q})\in \mathbb{R}_+^{p+q},$ then we
have
\begin{equation}
G^{\theta (\mathbf{\boldsymbol{\tau}}^{(p)})}(\mathbf{x}^{(p)})G^{\theta (\boldsymbol{\tau}^{(q)})}(\mathbf{x}^{(q)})\leq G^{\theta (\boldsymbol{\tau}^{(p+q)})}(\mathbf{x}^{(p+q)})\leq \min \{ G^{\theta
(\boldsymbol{\tau}^{(p)})}(\mathbf{x}^{(p)}),G^{\theta (\boldsymbol{\tau}^{(q)})}(\mathbf{x}^{(q)})\} ,  \label{eq2}
\end{equation}
for each $\mathbf{x}^{(p+q)}\in \mathbb{R}^{p+q}$ and $\tau_j=-\log F(x_j),$ $j=1,\ldots ,p+q.$

 The inequality on the right holds true for every
multivariate distribution, while the inequality on the left is a property
of MEV distributions. The lower bound corresponds to the case where
$\mathbf{Y}^{(p)}$ and $\mathbf{Y}^{(q)}$ are independent and the
upper bound corresponds to the case where $\mathbf{Y}^{(p)}$ and
$\mathbf{Y}^{(q)}$ are totally dependent.

From (\ref{eq2}) we obtain the following bounds for the multivariate extremal index function
$\theta(\boldsymbol{\tau}^{(p+q)}),\boldsymbol{\tau}^{(p+q)} \in
\mathbb{R}_+^{p+q}$.
\begin{equation}
\frac{\max \{ \theta (\boldsymbol{\tau}^{(p)})\gamma (\boldsymbol{\tau}
^{(p)}),\theta (\boldsymbol{\tau}^{(q)})\gamma (
\boldsymbol{\tau}^{(q)})\}}{\gamma (\boldsymbol{\tau}^{(p+q)})}\leq\theta({\boldsymbol{\tau}}^{(p+q)})\leq
\frac{\theta ({\boldsymbol{\tau}}^{(p)})\gamma ({\boldsymbol{\tau}}^{(p)})+\theta ({\boldsymbol{\tau}}^{(q)})\gamma (
{\boldsymbol{\tau}}^{(q)})}{\gamma ({\boldsymbol{\tau}}^{(p+q)})} , \label{eq3}
\end{equation}
where
\begin{equation*}
\gamma ({\boldsymbol{\tau}}^{(p+q)})=-\log G(F^{-1}(e^{-\tau_1}),\ldots ,F^{-1}(e^{-\tau_{p+q}}))=\lim_{n\rightarrow\infty}nP\left({\bf X}_1^{(p+q)} \not\leq u_n^{({\boldsymbol{\tau}}^{(p+q)})}\right),
\end{equation*}
\begin{equation*}
\gamma({\boldsymbol{\tau}}^{(p)})=\!\!\!\!\lim_{{\tau_j\rightarrow
0^+}\atop{j=p+1,\ldots ,p+q}} \gamma ({\boldsymbol{\tau}}^{(p+q)})\ \ \textrm{and}\ \
\gamma({\boldsymbol{\tau}}^{(q)})=
\lim_{{\tau_j\rightarrow 0^+}\atop{j=1,\ldots
,p}}\gamma ({\boldsymbol{\tau}}^{(p+q)}).
\end{equation*}

\noindent The next result follows naturally from these bounds.

\begin{proposition}
Suppose that $Q\in D(G)$ and $\mathbf{X}$ has
extremal index $\theta (\boldsymbol{\tau}^{(p+q)}),$
$\boldsymbol{\tau}^{(p+q)}\in \mathbb{R}_+^{p+q}.$

\noindent {\bf{(i)}} If $\ \widehat{\mathbf{Y}}^{(p)}$ and $\widehat{\mathbf{Y}}^{(q)}$  are independent, then $\mathbf{Y}^{(p)}$ and $\mathbf{Y}^{(q)}$ are independent if and only if
\begin{equation*}
\theta (\boldsymbol{\tau}^{(p+q)})=\frac{\theta (
\boldsymbol{\tau}^{(p)})\gamma (\boldsymbol{\tau}^{(p)})+\theta (
\boldsymbol{\tau}^{(q)})\gamma (\boldsymbol{\tau}^{(q)})}{
\gamma (\boldsymbol{\tau}^{(p)})+\gamma (\boldsymbol{\tau}^{(q)})},\ \boldsymbol{\tau}^{(p+q)}\in \mathbb{R}_+^{p+q}.
\end{equation*}

\noindent {\bf{(ii)}} If $\ \widehat{\mathbf{Y}}^{(p)}$ and
$\widehat{\mathbf{Y}}^{(q)}$ are totally dependent, then
$\mathbf{Y}^{(p)}$ and $\mathbf{Y}^{(q)}$ are totally dependent if
and only if $$\theta(\boldsymbol{\tau}^{(p+q)})=\frac{\max \{
\theta(\boldsymbol{\tau}^{(p)})\gamma(\boldsymbol{\tau}^{(p)}),
\theta(\boldsymbol{\tau}^{(q)})\gamma(\boldsymbol{\tau}^{(q)})\}
}{\max \{\gamma( \boldsymbol{\tau}^{(p)}),\gamma(\boldsymbol{\tau}^{(q)})\} } ,\ \boldsymbol{\tau}^{(p+q)}\in
\mathbb{R}^{p+q}.
$$
\end{proposition}

The necessary and sufficient conditions for $\mathbf{Y}^{(p)}$ and $\mathbf{Y}^{(q)}$
to be independent or totally dependent given in the previous result demand the evaluation of the extremal index function $\theta(\boldsymbol{\tau}^{(p+q)}),$ in each point $\boldsymbol{\tau}^{(p+q)}\in
\mathbb{R}_+^{p+q}.$ Nevertheless this task can be simplified with the characterizations, given in
Ferreira (\cite{Ferreira}), for independence and total dependence of the multivariate marginals of a
MEV distribution. These characterizations are essential to prove the following propositions which guarantee that the independence or total dependence between $\mathbf{Y}^{(p)}$ and $\mathbf{Y}^{(q)}$ only depends on the value of the extremal
index in some points.

\begin{proposition}
Suppose that $Q\in D(G)$ and the sequence $\mathbf{X}=\{\mathbf{X}_{n}^{(p+q)}\} _{n\geq 1}$ has extremal index $\theta(\boldsymbol{\tau}^{(p+q)}),$ $\boldsymbol{\tau}^{(p+q)}\in
\mathbb{R}_+^{p+q}.$

The sub-vectors $\mathbf{Y}^{(p)}$ and $\mathbf{Y}^{(q)}$ are independent if and only if
\begin{equation}
\label{eq4}
\theta (\mathbf{1}^{(p+q)})=\frac{\theta (\mathbf{1}^{(p)})\gamma (\mathbf{1}^{(p)})+
\theta (\mathbf{1}^{(q)})\gamma (\mathbf{1}^{(q)})}{\gamma (\mathbf{1}^{(p+q)})},
\end{equation}
where $\mathbf{1}^{(k)}=(1,\ldots,1),\ k>1$, denotes the
$k$-dimensional unitary vector.
\end{proposition}

\bdem
Suppose that $\mathbf{Y}^{(p)}$ and $\mathbf{Y}^{(q)}$ are independent. Since (\ref{eq3}) holds for all
\linebreak $\mathbf{\boldsymbol{\tau}}^{(p+q)}\in \mathbb{R}_+^{p+q}$, we have
in particular for $ \boldsymbol{\tau}^{(p+q)}=(\tau,\ldots ,\tau)\in \mathbb{R}_+^{p+q}$, with \linebreak $\tau \equiv \tau(x)=-\log F(x),\ x
\in \mathbb{R},$
$$
\theta(\boldsymbol{\tau}^{(p+q)})=\frac{\theta(\boldsymbol{\tau}^{(p)})\gamma(\boldsymbol{\tau}^{(p)})+\theta(\boldsymbol{\tau}^{(q)})\gamma (\boldsymbol{\tau}^{(q)})}{\gamma(\boldsymbol{\tau}^{(p+q)})}.
$$
\noindent Now attending to the fact that $\theta (c\boldsymbol{\tau}^{(k)})=\theta
(\boldsymbol{\tau}^{(k)})$ for each $\boldsymbol{\tau}^{(k)}\in
\mathbb{R}_+^k,$ $k>1$ and $c>0$, we can write
$$
\theta(\boldsymbol{\tau}^{(p+q)})=\theta (\mathbf{1}^{(p+q)}),\quad \theta(\boldsymbol{\tau}^{(p)})=\theta(\mathbf{1}^{(p)}),\quad \theta
(\boldsymbol{\tau}^{(q)})=\theta(\mathbf{1}^{(q)}),
$$
and  from (\ref{eq1}), for all
$\boldsymbol{\tau}^{(p+q)}=(\tau ,\ldots ,\tau )\in
\mathbb{R}_+^{p+q}$,
\begin{eqnarray}
\gamma({\boldsymbol{\tau}}^{(p+q)}) &=&-\log G(F^{-1}(e^{-\tau }),\ldots ,F^{-1}(e^{-\tau }))=-\log D_{G}(e^{-\tau},\ldots ,e^{-\tau}) \nonumber \\
&=&-\log D_G^{\tau}(e^{-1},\ldots ,e^{-1})=\tau \gamma (\mathbf{1}^{(p+q)}),
\label{eq5}
\end{eqnarray}
$\gamma(\boldsymbol{\tau}^{(p)})=\tau \gamma
(\mathbf{1}^{(p)})\text{ and  }\gamma(\boldsymbol{\tau}^{(q)})=\tau
\gamma (\mathbf{1}^{(q)}).$ Equality (\ref{eq4}) is now
straightforward.

On the other hand  if (\ref{eq4}) is verified, then for
$\mathbf{x}^{(p+q)}=(x,\ldots ,x)$ we have
\begin{eqnarray*}
G_{\mathbf{Y}}(\mathbf{x}^{(p+q)}) &=&G^{\theta (\mathbf{1}^{(p+q)})}(\mathbf{x}^{(p+q)})=D_{G}^{\theta (\mathbf{1}^{(p+q)})}(e^{-\tau},\ldots
,e^{-\tau})\\
&=&D_{G}^{\theta(\mathbf{1}^{(p+q)})\tau}(e^{-1},\ldots
,e^{-1}) =\exp(-\tau \gamma (\mathbf{1}^{(p+q)})\theta (\mathbf{1}^{(p+q)}))\\
&=&\exp (-\tau(\theta(\mathbf{1}^{(p)})\gamma (\mathbf{1}^{(p)})+\theta
(\mathbf{1}^{(q)})\gamma(\mathbf{1}^{(q)})))=G_{\mathbf{Y}^{(p)}}(\mathbf{x}^{(p)})G_{\mathbf{Y}^{(q)}}(\mathbf{x}^{(q)})
\end{eqnarray*}
and from Proposition 2.1 (\cite{Ferreira}) we conclude that $\mathbf{Y}^{(p)}$ and
$\mathbf{Y}^{(q)}$ are independent.\edem

\begin{proposition}
Suppose that $Q\in D(G)$, $\mathbf{X}=\{\mathbf{X}_{n}^{(p+q)}\} _{n\geq 1}$ has extremal index\linebreak $\theta(
\boldsymbol{\tau}^{(p+q)}),$ $\boldsymbol{\tau}^{(p+q)}\in \mathbb{R}_+^{p+q}.$

\noindent {\bf{(i)}} If $\mathbf{Y}^{(p)}$ and $\mathbf{Y}^{(q)}$ are totally dependent then
there exists  $\boldsymbol{\tau}^{(p+q)}\in \mathbb{R}_+^{p+q}$ with \linebreak
$\tau_j\equiv \tau_j(x_j)=-\log F(x_j), x_j \in \mathbb{R},\  j=1,\ldots,p+q$, such that
$$\gamma(\boldsymbol{\tau}^{(p)})\theta(\boldsymbol{\tau}^{(p)})=\gamma(\boldsymbol{\tau}^{(q)})\theta(
\boldsymbol{\tau}^{(q)})=\theta_1\tau_1\ldots=\theta_{p+q}\tau_{p+q}=d>0$$
and $\theta(\boldsymbol{\tau}^{(p+q)})=\gamma\left(\frac{\boldsymbol{\tau}^{(p+q)}}{d}\right)^{-1}$.

\noindent {\bf{(ii)}} If there exists $\boldsymbol{\tau}^{(p+q)}\in \mathbb{R}_+^{p+q}$ with
$\tau_j\equiv \tau_j(x_j)=-\log F(x_j),\ x_j \in \mathbb{R},\  j=1,...,p+q$, such that
\begin{equation*}
\gamma (\boldsymbol{\tau}^{(p+q)})\theta (\boldsymbol{\tau}^{(p+q)})=\theta_1\tau_1\ldots=\theta_{p+q}\tau_{p+q}=d>0,
\end{equation*}
then $\mathbf{Y}^{(p)}$ and $\mathbf{Y}^{(q)}$ are totally dependent.
\end{proposition}

\bdem
{\bf{(i)}} From Proposition 2.1 (\cite{Ferreira}), if $\mathbf{Y}^{(p)}$ and
$\mathbf{Y}^{(q)}$ are totally dependent, then there exists
$\mathbf{\boldsymbol{\tau}}^{(p+q)}\in \mathbb{R}_+^{p+q}$ such that
$$
\theta(\boldsymbol{\tau}^{(p)})\gamma(\boldsymbol{\tau}^{(p)})=\theta(\boldsymbol{\tau}^{(q)})\gamma(
\boldsymbol{\tau}^{(q)})=\theta (\boldsymbol{\tau}^{(p+q)})\gamma(\boldsymbol{\tau}^{(p+q)})=d=\theta_1\tau_1=\ldots=\theta_{p+q}\tau_{p+q},
$$
with $d \in ]0,1[.$ Hence
\begin{eqnarray*}
\theta(\boldsymbol{\tau}^{(p+q)}) &=&\frac{d}{\gamma(\boldsymbol{\tau}^{(p+q)})}=
\frac{d}{-\log D_G(\exp (-\tau_1),\ldots ,\exp (-\tau_{p+q}))}\\
&=&\frac{1}{-\log D_{G}\left( \exp \left( -\frac{\tau_{1}}{d}\right),\ldots ,\exp \left(-\frac{\tau_{p+q} }{d}\right)\right)}=\frac{1}{\gamma\left(\frac{\boldsymbol{\tau}^{(p+q)}}{d}\right)}.
\end{eqnarray*}
\edem

Another way to look at issues concerning independence or total dependence is to use parameters that measure the strength of dependence between ${\bf{Y}}^{(p)}$ and ${\bf{Y}}^{(q)}$. We therefore define, in the following result, the dependence structure of ${\bf{Y}}^{(p)}$ and ${\bf{Y}}^{(q)}$ through the coefficient $\epsilon^{({\bf{Y}}^{(p)},{\bf{Y}}^{(q)})}$ of Ferreira (\cite{Ferreira}). This coefficient emerged from the extremal coefficient of $\bf Y$, $\epsilon^{\bf{Y}}$, defined in Martins and Ferreira (\cite{Martins}) as
$$G^{\theta(\textbf{1}^{(p+q)})}(\textbf{x}^{(p+q)})=F^{{\epsilon}^{\mathbf{Y}}}(x), \ \ x\in \mathbb{R},$$
and the relationship
$$P\left({\bf{Y}}^{(p)}\leq{\textbf{x}}^{(p)},{\bf{Y}}^{(q)}\leq{\textbf{x}}^{(q)}\right)=\left(G_{\mathbf{Y}}^{(p)}({\textbf{x}}^{(p)})
G_{\mathbf{Y}}^{(q)}({\textbf{x}}^{(q)})\right)^{\frac{\epsilon^{\mathbf{Y}}}{\epsilon^{{\bf{Y}}^{(p)}}+\epsilon^{{\bf{Y}}^{(q)}}}}.$$
It is then defined as
$\epsilon^{\left({\bf{Y}}^{(p)},{\bf{Y}}^{(q)}\right)}=\frac{\epsilon^{\mathbf{Y}}}{\epsilon^{{\bf{Y}}^{(p)}}+\epsilon^{{\bf{Y}}^{(q)}}}$ and has the following interesting pro\-per\-ties.

\begin{proposition}$\ $\\

\noindent {\bf{(i)}} $\epsilon^{({\bf{Y}}^{(p)},{\bf{Y}}^{(q)})}=\frac{\theta({\bf{1}}
^{(p+q)})\gamma ({\bf{1}}^{(p+q)})}{\theta({\bf{1}}^{(p)})\gamma(
{\bf{1}}^{(p)})+\theta({\bf{1}}^{(q)})\gamma({\bf{1}}^{(q)})}$

\noindent {\bf{(ii)}} $\epsilon^{({\bf{Y}}^{(p)},{\bf{Y}}^{(q)})}=1$ if and only if ${\bf{Y}}^{(p)}$
and ${\bf{Y}}^{(q)}$ are independent.

\noindent {\bf{(iii)}} If $\ {\bf{Y}}^{(p)}$ and ${\bf{Y}}^{(q)}$  are totally dependent, then
$\epsilon^{({\bf{Y}}^{(p)},{\bf{Y}}^{(q)})}= \frac{\max\{\epsilon^{{\bf{Y}}^{(p)}},\epsilon^{{\bf{Y}}^{(q)}}\}}{\epsilon^{{\bf{Y}}^{(p)}}+\epsilon^{{\bf{Y}}^{(q)}}}.$

\end{proposition}

\bdem
\noindent {\bf{(i)}} Since
$$
\epsilon^{({\bf{Y}}^{(p)},{\bf{Y}}^{(q)})}=\frac{\theta({\bf{1}}
^{(p+q)})\log G({\bf x}^{(p+q)})}{\theta({\bf{1}}^{(p)}) \log G({\bf x}^{(p)})+\theta({\bf{1}}^{(q)})\log G({\bf x}^{(q)}))}
$$
the result follows from (\ref{eq5}).

\noindent {\bf{(ii)}} It is an immediate consequence of (i) and Proposition 2.

\noindent{\bf{(iii)}} If ${\bf{Y}}^{(p)}$ and ${\bf{Y}}^{(q)}$ are totally dependent then from (\ref{eq2}), we have
\begin{equation*}
\epsilon^{({\bf{Y}}^{(p)},{\bf{Y}}^{(q)})} = \frac{-\log \min \{G^{\theta (\bf{1}^{(p)} )}({\bf{x}}^{(p)}),G^{\theta ({\bf{1}}^{(q)})}({\bf{x}}^{(q)})\}}
{\theta({\bf{1}}^{(p)}) \log G({\bf{x}}^{(p)})+ \theta({\bf{1}}^{(q)}) \log G({\bf{x}}^{(q)})}=\frac{\max\{\epsilon^{{\bf{Y}}^{(p)}},\epsilon^{{\bf{Y}}^{(q)}}\}}{\epsilon^{{\bf{Y}}^{(p)}}+\epsilon^{{\bf{Y}}^{(q)}}}.
\end{equation*}
\edem

\section{Examples}
\label{sec:7}

\begin{example}
{\rm{\noindent Let $\{ Y_n \} _{n\geq 1}$ be a sequence of i.i.d.
random variables with common d.f. $F$ and consider an
auto-regressive sequence of maxima $ \{X_n\}_{n\geq 1}$ defined by
$$X_n=\max\{Y_n,Y_{n+1}\},\ n\geq 1,$$
with marginal distribution function $F^2.$

Let $\{u_n^{(\tau_i)}\} _{n\geq 1},\ i=1,\ldots,p,$ and $\{v_n^{(\tau'_j)}\} _{n\geq 1},j=p+1,\ldots ,p+q,$
be sequences of real numbers such that $\lim_{n\to \infty}n(1-F^{2}(u_{n}^{(\tau _{i})}))=\tau _{i}$ and $\lim_{n\to \infty}nF^2(-v_n^{(\tau'_j)})) =\tau_{j}^{\prime}.$

The sequences $\{X_n\}_{n\geq 1}$ and $\{-X_n\}_{n\geq 1}$ have, respectively, extremal indexes $\theta_1=1/2$ and $\theta_2=1.$

\noindent For sequences $\mathbf{X}_{n}^{(p+q)}=\left\{
\begin{array}{lcl}
X_{n,i}=X_{n}&,& i=1,\ldots p\\
X_{n,i}=-X_{n}&,& i=p+1,\ldots ,p+q
\end{array}
\right.$, $\mathbf{X}_{n}^{(p)}=(X_{n},\ldots ,X_{n})$ and $\mathbf{X}_n^{(q)}=(-X_n,\ldots ,-X_n)$, we have
\begin{eqnarray*}
\lim_{n\to \infty}P(\mathbf{M}_{n}^{(p)}\leq \mathbf{u}_n^{(\boldsymbol{\tau}^{(p)})})&=&\exp
\left(-\frac{1}{2}\underset{1\leq j\leq p}{\max}\tau_{j}\right),\\
\lim_{n\to \infty}P(\widehat{\mathbf{M}}_{n}^{(p)}\leq \mathbf{u}_n^{(\boldsymbol{\tau}^{(p)})})&=&\exp \left( -\underset{1\leq j\leq p}{
\max }\tau _{j}\right) , \\
\lim_{n\to \infty}P\left( \mathbf{M}_{n}^{(q)}\leq (v_{n}^{(\tau _{p+1}^{\prime })},\ldots
,v_{n}^{(\tau _{p+q}^{\prime })})\right) &=&\exp \left( -\underset{
p+1\leq j\leq p+q}{\max }\tau _{j}^{\prime }\right) .
\end{eqnarray*}
Since the order statistics maximum and minimum are asymptotically
independent (\cite{Davis},\cite{Pereira}) we obtain
\begin{equation*}
{P( \mathbf{M}_n\leq(u_n^{(\tau_1)},\ldots ,u_n^{(\tau
_p)},v_n^{(\tau_{p+1}^{\prime })},\ldots ,v_n^{(\tau _{p+q}^{\prime
})}) ) }\xrightarrow [n\rain]{}\exp \left( -\frac{1}{2}\underset{1\leq
j\leq p}{\max }\tau _{j}-\underset{p+1\leq j\leq
p+q}{\max }\tau _{j}^{\prime }\right)
\end{equation*}
and consequently  $\gamma(\boldsymbol{\tau}^{(p+q)})=\gamma(\boldsymbol{\tau}^{(p)})+\gamma(\boldsymbol{\tau}^{(q)})=\underset{1\leq j\leq p}{\max }\tau _{j}+\underset{
p+1\leq j\leq p+q}{\max }\tau _{j}^{\prime }$   and   \linebreak
$\theta(\boldsymbol{\tau}^{(p+q)})\gamma(\boldsymbol{\tau}^{(p+q)}) =\frac{1}{2}\ \underset{1\leq j\leq p}{\max }\tau _{j}+\underset{p+1\leq j\leq p+q}{\max}
\tau_{j}^{\prime }.$ Therefore
\begin{equation*}
\theta (\boldsymbol{\tau}^{(p+q)})=\frac{\theta (\boldsymbol{\tau
}^{(p)})\gamma ( \boldsymbol{\tau}^{(p)})+\theta (\boldsymbol{\tau
}^{(q)})\gamma (\boldsymbol{\tau} ^{(q)})}{\gamma (\boldsymbol{\tau
}^{(p)})+\gamma (\boldsymbol{\tau}^{(q)})}.
\end{equation*} }}
\end{example}

\begin{example}
{\rm{\noindent Let ${\bf{U}}=\{U_n\}_{n\geq 1}$ be a sequence of i.i.d. random variables with common d.f. $H$ in the domain of attraction of the extreme value distribution $F$, and independent of the i.i.d. chain ${\bf{J}}=\{J_n\}_{n\geq 1}$ such that $P(J_1=0)=P(J_1=1)=1/2$.

Let us consider a stationary 1-dependent sequence ${\bf{Z}}=\{Z_n\}_{n\geq 1},$ defined as $Z_n=U_n$ if $J_n=0$ and $Z_n=U_{n+1}$ otherwise, and let ${\bf{v}}=\{v_n\}_{n\geq1}$ be a sequence of normalized levels to $\bf{Z}$, and consequently also to $\bf{U}$.

We can now define a 3-dependent stationary sequence
${\bf{X}}=\{{\bf{X}}_n=(X_{n,1},X_{n,2},X_{n,3})\}$ as
$$(X_{n,1},X_{n,2},X_{n,3})=(Z_n,Z_{n+2},Z_{n+1}), \ \ n\geq 1,$$
with common distribution function
$$T(x_1,x_2,x_3)=\frac{1}{2}\prod_{i=1}^3H(x_i)+\frac{1}{4}H(x_1)H(\min\left\{x_2,x_3\right\})+\frac{1}{4}H(x_2)H(\min\left\{x_1,x_3\right\})$$
belonging to the domain of attraction of
$$G(x_1,x_2,x_3)=\left\{
\begin{array}{lcl}
F(x_1)F(x_2)F^{\frac{1}{2}}(x_3) & , & \ x_1<x_3\wedge x_2<x_3 \\
F(x_1)F^{\frac{3}{4}}(x_2)F^{\frac{3}{4}}(x_3) & , &  \ x_1<x_3\wedge x_3\leq x_2 \\
F^{\frac{3}{4}}(x_1)F(x_2)F^{\frac{3}{4}}(x_3) & , &  \ x_3\leq x_1\wedge x_2< x_3 \\
F^{\frac{3}{4}}(x_1)F^{\frac{3}{4}}(x_2)F(x_3) & , &  \ x_3\leq x_1\wedge x_3 \leq x_2
\end{array}
\right.$$

Now applying Proposition 2.1 of Smith and Weissman (1996) (\cite{smith}) to the sequence
$U^{\bf{X}}=\{\max\{X_{n,1},X_{n,2},X_{n,3}\}=\max\{Z_n,Z_{n+1},Z_{n+2}\}\}_{n\geq1}$ which verifies the condition $D^{(k)}(v_n), \ k=2$, of Chernick et al. (\cite{Chernick}), we easily obtain
$$\theta^{\bf{X}}(1,1,1)=\lim_{n\rightarrow \infty}\frac{P\left(\max\left\{Z_1,Z_2,Z_3\right\}>v_n\geq \max \left\{Z_2,Z_3,Z_4\right\}\right)}{P\left(\max\left\{Z_1,Z_2,Z_3\right\}>v_n\right)}=\frac{3}{10}.$$

 For random vectors $\widehat{{\bf{Y}}}=\left(\widehat{Y}_1,\widehat{Y}_2,\widehat{Y}_3\right)$ and ${\bf{Y}}=\left(Y_1,Y_2,Y_3\right)$  with d.f. $G_{\widehat{{\bf{Y}}}}\equiv G$ and $G_{\bf{Y}}\equiv G^{\theta (\boldsymbol{\tau}^{(3)})}$, ${\bf{Y}}^{(2)}=\left(Y_1,Y_2\right)$ and ${\bf{Y}}^{(1)}=Y_3$ we obtain $\epsilon^{\widehat{{\bf{Y}}}}=\frac{5}{2}$, $\epsilon^{\bf{Y}}=\frac{3}{4}$, $\epsilon^{{\widehat{\bf{Y}}}^{(2)}}=2$, $\epsilon^{{\widehat{{\bf{Y}}}}^{(1)}}=1$, $\epsilon^{{{\bf{Y}}}^{(1)}}=\frac{3}{4}=\epsilon^{{{\bf{Y}}}^{(2)}}$. Consequently  $\epsilon^{\left({\bf{Y}}^{(2)},{\bf{Y}}^{(1)}\right)}=\frac{1}{2},$ $\epsilon^{\left(\widehat{\bf{Y}}^{(2)},\widehat{\bf{Y}}^{(1)}\right)}=\frac{5}{6}$ and from Proposition 4 we conclude that neither ${\bf{Y}}^{(1)}$ and ${\bf{Y}}^{(2)}$ nor $\widehat{\bf{Y}}^{(1)}$ and $\widehat{\bf{Y}}^{(2)}$ are independent.

Nevertheless, there exists $\boldsymbol{\tau}^{(3)}=(1,1,1)$ such that $\gamma(\boldsymbol{\tau}^{(3)})\theta(\boldsymbol{\tau}^{(3)})=\theta_1\tau_1=\theta_2\tau_2=\theta_3\tau_3=\frac{3}{4}$ and attending to Proposition 3 we can say that ${\bf{Y}}^{(1)}$ and ${\bf{Y}}^{(2)}$ are totally dependent.}}
\end{example}

%
%

\end{document}